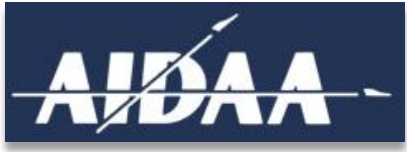

# MULTIDISCIPLINARY OPTIMIZATION FOR GAS TURBINES DESIGN

F. Bertini[1], L. Dal Mas[2], L. Vassio[3] and E. Ampellio[3*]

[1]AVIO S.p.A., R&TD, 10040 – Rivalta di Torino (TO), Italy

[2]Università di Padova, Dipartimento di Ingegneria Industriale, 35100 – Padova (PD), Italy

[3]Politecnico di Torino, Dipartimento di Ingegneria Meccanica e Aerospaziale, 10129 – Torino (TO), Italy

*enrico.ampellio@polito.it

**ABSTRACT**

*State-of-the-art aeronautic Low Pressure gas Turbines (LPTs) are already characterized by high quality standards, thus they offer very narrow margins of improvement. Typical design process starts with a Concept Design (CD) phase, defined using mean-line 1D and other low-order tools, and evolves through a Preliminary Design (PD) phase, which allows the geometric definition in details. In this framework, multidisciplinary optimization is the only way to properly handle the complicated peculiarities of the design.*

*The authors present different strategies and algorithms that have been implemented exploiting the PD phase as a real-like design benchmark to illustrate results. The purpose of this work is to describe the optimization techniques, their settings and how to implement them effectively in a multidisciplinary environment. Starting from a basic gradient method and a semi-random second order method, the authors have introduced an Artificial Bee Colony-like optimizer, a multi-objective Genetic Diversity Evolutionary Algorithm [1] and a multi-objective response surface approach based on Artificial Neural Network, parallelizing and customizing them for the gas turbine study. Moreover, speedup and improvement arrangements are embedded in different hybrid strategies with the aim at finding the best solutions for different kind of problems that arise in this field.*

## 1   INTRODUCTION

Recent international regulations like ACARE [2] for Europe and ICAO [3] worldwide impose to reduce fuel consumption, pollutant and noise emissions for greener aircraft and engines. Focusing on propulsion system, the intensive application and tuning of high-performance and multi-objective optimization strategies is imperative. Especially as regards Ultra High Bypass Turbofan (UHBT), the design of Low Pressure Turbine (LPT) is a critical and challenging task, since this component has a great impact on Specific Fuel Consumption and only a multi-disciplinary approach can allow obtaining the best configuration.

Low pressure turbine module is a very complex system, made of several physical components and operating in a fluid environment difficult to be numerically simulated in its completeness. Therefore, the design of a LPT is an intrinsically multidisciplinary problem, normally characterized by a high dimensionality in terms of both free-to-change parameters, and objective functions. Since the early CD phase, there are many important characteristics that need to be evaluated at the same time (e.g. aerodynamic performances, structural properties, acoustic behaviour and different heterogeneous constraints). The correct definition of the optimization problem is very important, given that the costs of changes in a project are higher when advancing in the design process and available degrees of freedom are fewer.

In general, input/output definitions and computational tools become more and more detailed and specific moving from CD to Detailed Design (DD). For example, in concept design, only fundamental LPT traits are considered, like flow path, blade number and work split in order to respect mechanical and aero-acoustic requirements through mean-line 1D correlation-based solvers. On the other hand, in DD some deeply specialized characteristics, like local geometries of tip clearance or non-axisymmetric endwall contouring are optimized in order to essentially refine the overall performance, using 3D full Navier-Stokes multistage unsteady solvers.

From these simple considerations, it is apparent that many kind of optimization problems can arise during the design of a LPT. Quasi mono-objective problems with huge and jagged domain of parameters need to be addressed, as well as inherently multi-objective ones with strongly contrasting targets, but easy-knowable boundaries to set for feasibility. In this scenario, a very large spectrum of optimization algorithms and methodologies, each having its potentials and limitations, can be applied. While some of them are general and easy to be used, others need to be properly tuned in order to successfully complete the optimization. Therefore, the perfect knowledge of the physical/engineering problem to be solved is required to choose the best algorithm or methodology.

To emphasize what has been said so far, the authors will introduce and analyse different optimization strategies, from the standard to the forefront ones. Aeronautical LPTs optimization in PD phase represents a real-like test case for meaningful comparisons. The final aim of the paper is to provide a straightforward compendium of modern cutting edge techniques easily extendable and applicable to many research areas.

## 2    DESIGN PROCESS

This work fits in the standard procedure for aero design of aeronautic Low Pressure Turbines currently adopted by AVIO S.p.A.. Figure 1 reports the complete framework which includes both the CD and the PD phases. Each block of the flowchart was thoroughly described in [4].

First of all, the Concept Design phase is basically carried out using a 1D Meanline Multidisciplinary solver, embedded in the in-house tool called 'Tùrbine' together with various optimization methods herein discussed. This tool is able to manage the overall turbine layout by defining some fundamental parameters like operating thermodynamic cycle, encumbrance geometries and flow path. 'Tùrbine' calculates velocity triangles by solving Euler equations and it exploits revised classical loss correlations for performance estimation [5], finally building the preliminary turbine cross section. Within 'Tùrbine', it is also possible to carried out mono and multi–objective optimizations by properly setting the design variables and objectives.

2D solution represents the second step in CD. This part of the design is aimed at optimizing the preliminary spanwise work and outlet aerodynamic angle distribution, while meanline values are maintained constant. The through-flow computations follow the same approach of meanline ones, except that many streamline are considered along the radial extension of the flow path.

Blade shape is realized through the identification of a discrete number of sections, whose airfoils are parametrically defined. The 3D geometry of each blade can be automatically generated using the data from 1D meanline or 2D spanwise solution, or even including the benefits obtained from Q3D optimizations during PD phase described below. Database corrections, inverse correlation and free vortex principles are applied.

During the PD phase, the geometric refinement of 3D airfoils is achieved with the help of multidisciplinary optimization. Since employing directly three-dimensional simulations would be expensive in terms of both computational resources and time, few reference sections are usually studied by means of faster Q3D analyses [12].

In order to interpret the physical 3D phenomena that got lost in the numerical model simplification (e.g. secondary vortex effects), the flow channel diffusion has to be adjusted according to the 3D results. This operation is normally referred to as "fitting". Multi-row steady and unsteady calculation are available to rearrange velocity triangles and stator/rotor interaction using CFD; single row simulation are introduced for perfecting airfoils geometries. In both cases, at least aerodynamic and structural requirements have always to be satisfied.

CFD 3D multistage complete calculations [13] are performed after every optimization cycle to confirm the improvement made and obtain fully trustful reference results.

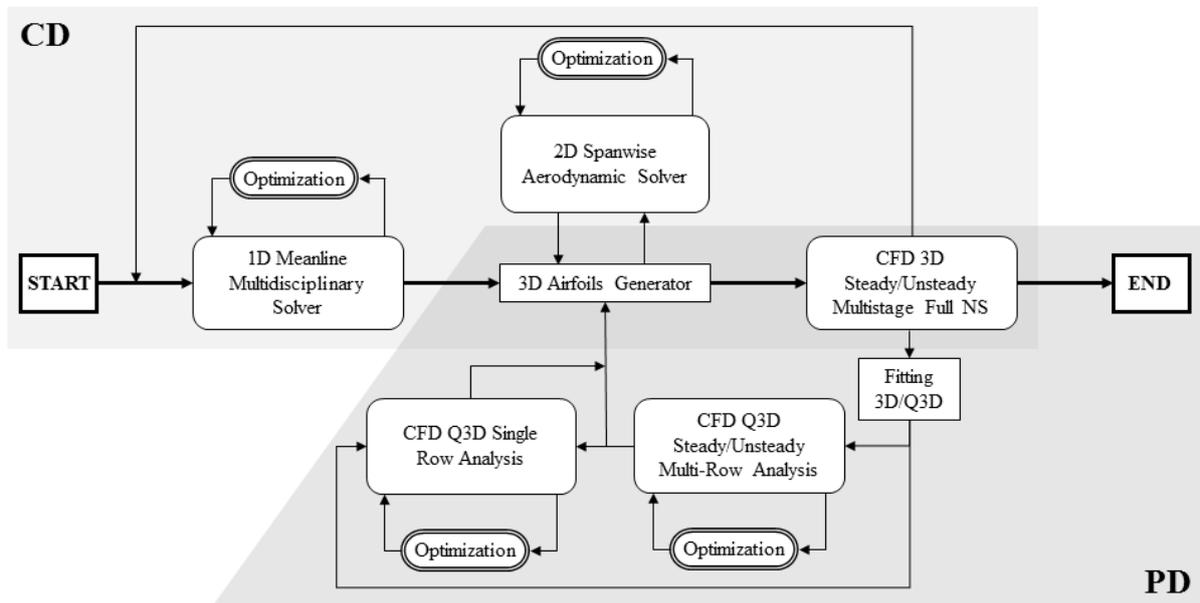

Figure 1: Typical AVIO S.p.A. aereo-design process for LPTs

## 3 OPTIMIZATION ALGORITHMS

As introduced in the first two chapters, various optimization methods are necessary in the different phases of the LPT design. In this section some valuable algorithms, used for comparisons in the next chapter, will be briefly outlined. They are presented, without loss of generality, as possible solution of a minimization problem. These algorithms were all implemented by the authors, with the only exception of the GeDEA, which was developed by the University of Padova.

### 3.1 Gradient Descent (GD)

Gradient Descent (GD or Steepest Descent) is a classical first-order iterative optimization algorithm for continuous differentiable functions [6]. It often represents a good first try for optimization methods thanks to its simplicity and velocity. To proceed towards the minimum from the current point, the method moves in the opposite direction of the gradient evaluated in the point. In the case where an analytical characterization of the gradient is not feasible, an estimation of it should be computed. In this work, in order to reduce the number of function evaluations, is used a forward finite difference method with first order accuracy, performed in parallel with as many evaluations as are the free parameters. Once the direction of improvement is found, the step size has to be accurately chosen to reach quickly a local minimum. In the studied case a line optimization search over the step size, implemented in parallel, is performed.

A negative aspect of GD is that there is no space investigation, following greedily the steepest descent to a close local minimum. Moreover, if the functions studied are noisy under a certain step threshold, the estimation of derivatives is hard to calibrate. Notice also that this algorithm is totally deterministic, not depending from random components, and automatically stops itself when no improvement direction is found.

### 3.2 Semi-Random Walk with second order Interpolation (IRW)

The algorithm proposed is an improvement of a simple random walk [6] over the parameters space. Iteratively from the best point found, a random step is performed and, if it is not an improvement, the opposite step is tried. If even such step is not an improvement then a simple estimation of the directional curvature of the function is computed through a second order interpolation of the three previous points. Then a new point, minimum of such parabola, is tested.

With this procedure the random walk can exit more easily from local minima with respect to the GD. Moreover, a particular experience on the problem is not required to implement the procedure. On the contrary this algorithm is very slow due to the totally random component that avoid to use data history smartly to find a direction of improvement. Additionally, it cannot be effectively parallelized due to the serial nature of the walk.

### 3.3 Artificial super-Bee enhanced Colony (AsBeC)

Artificial Bee Colony (ABC) is one of the most recent and promising evolutionary technique, inspired by the intelligent behaviour of honey bees. Originally introduced in 2005 by [7], it is now widely spread and studied in many research fields for its simplicity, effectiveness and robust behaviour. This algorithm embraces the principles of Particle Swarm Optimization (PSO) and Differential Evolution (DE). ABC provides a population-based search modus operandi in which bees fly around a multidimensional space. A group of bees searches food in the space randomly, while another one moves depending on the experience of all the nest mates. The food sources with higher nectar are constantly updated, memorized and exploited to look around, while poorer ones are instead abandoned. Hence ABC system combines local and global search methods. Several types of upgrading have already been proposed over the years to improve the speed, the quality and/or the explorative ability of the swarm. An example of current use in turbomachinery is reported in [8].

The authors developed and used an improved version of the original algorithm, called Artificial super-Bee enhanced Colony (AsBeC), to which will be soon dedicated an in-depth technical article. In short, AsBeC is a parallelized, hybrid and enhanced scheme. The bee movement here follows the aforementioned semi-random second order method (super-Bee principle). Bees repositioning is based on both nectar amount evolution and space filling properties. Moreover analyses of equal configurations are avoided. Finally the objective function evaluations are executed in parallel at each swarm movement iteration.

AsBeC is easy to handle thanks to the small number of tuneable characteristic parameters. As in the IRW, it still wastes many time expensive simulations, due to its implicit random nature.

## 3.4 Genetic algorithm (GeDEA)

The genetic algorithms (GAs) are a particular class of evolutionary algorithms inspired by the natural evolutions. Nowadays they are widely used for the solution of multi-objective optimization problems, since they tend to be good at finding generally good global solutions. The Genetic Diversity Evolutionary Algorithm (GeDEA) [1] was developed at the University of Padova (UNIPD) and it was implemented within the AVIO procedures for the optimization of LPTs, under a collaborative project between the university and the company.

The GeDEA starts creating a random population of a given number of individuals. A mating pool is created form the initial population, then the offspring is generated by crossover and mutation. The presence of clones is automatically avoided by the algorithm, which substitutes the copies with randomly generated individuals, thus encouraging the exploration of the search space. The new individuals are processed to evaluate the objective functions. The algorithm estimates the genetic diversity of each individual of the current population with respect to the other, by means of a common distance metric. Then, the parents and the offspring are ranked using the non-dominating sorting procedure. Therefore, the GeDEA provides a dual selection pressure towards both the Pareto optimal set and the maintenance of the diversity within the population. According to the assigned ranks, the best solutions are selected for survival and for the generation of the successive mating pool, while the other individuals are suppressed. It is worth noting that the original version of the algorithm was partially modified to guarantee the survival of those individual characterized by the highest weighted sum of the selected objective functions (elitism). The optimization loop continues until the number of generations prescribed by the user is reached.

At the end of the multi-objective optimization process, the algorithm provides a set of optimal solutions, according to the Pareto concept, among which the user can choose the one that is more suitable for his needs. The algorithm is easily parallelizable, since each individual belonging to a specific generation can be analysed separately from the other. Another advantage of GAs is the widespread exploration of the feasible design space, the bounds of which should be accurately defined in order for the algorithm to work properly and efficiently.

## 3.5 Artificial Neural Network (ANN) applied to Latin Optimal Hypercube (LOH) sampling

Artificial Neural Networks are mathematical models inspired by biological neural networks and are used to model complex relationships between inputs and outputs [9]. An ANN is formed by an interconnected cluster of artificial computational units, called neurons, that process information through signal alteration and recombination. Once the architecture of the network is set, the model adapts itself during a training phase. In the studied case, the ANN is trained by the standard back-propagation gradient-based procedure [9]. Thanks to the training phase an ANN can estimate outputs corresponding to new inputs never tried before, acting as a 'black-box'. The incredible speed in outcomes computed through the network is one of its strengths, with respect to the direct simulation of the emulated phenomenon.

The scheme developed by the authors is intended to use the ANN as a fast surrogate method for the estimation of the objective functions, without using direct Q3D analyses. Therefore, ANN approach permits to manage multi-objective optimizations and leads easily to an acceptable approximation of the Pareto Front. As in the GeDEA algorithm, the choice of many parameters still remain hard and a certain experience on initial space boundaries is needed.

The key point for obtaining good results is the selection of the training points: they are taken from a well-distributed sampling, the LOH [10], within the investigation ranges. The dimensions of the hypercube has been chosen to suitably address the optimization problem with the aim at maximizing quality to computation resources ratio. However, high dimensionality in parameters would be a limitation for the velocity of this scheme. Best Pareto points suggested by the network are then directly simulated and used to perform subsequent refinement cycles of LOHs+ANNs, reducing investigation ranges nearby lowest weighted sum point.

### 3.6 Accelerating techniques

Some important accelerating technologies have been implemented within the aforementioned optimization algorithms. In addition to probabilistic choice in parameter modification and range resize [11], the authors developed a prediction routine that works as a trend analyser. This routine exploits the knowledge about the history of the best solutions found and how the variables changed: given the last best solution values and its relative parameters, the algorithm guesses a new point to try, computing a weighted average of the last directions of improvement. Such weights depends both on the quality of the improvements and on the proximity to the current solution. This accelerator works better with the algorithms related to a path in the search space (IRW, AsBeC).

## 4 TEST CASE: HIGH SPEED TURBINE

The real-like test case selected for comparison purposes is a rotor blade of a high-speed turbine configuration. This LPT layout is typically suited for geared-fan or open rotor engines and it requires a strong and careful aero-mechanic optimization procedure integrated in the Q3D analysis tools. While mechanical strength is predominant at the hub section, aerodynamic features are critical at the tip section. This background outlines differential optimization problems depending on the blade span. The distinctiveness traits just illustrated cause the 3D blade shape constructed by the profile generator after the first CD (see Figure 1 and Part 2) to be inadequate for the case under study. For this reason fast and effective multi-objective optimizations are fundamental.

Keeping in mind the above, the authors will present one complete design cycle in PD phase, consisting in the 3D/Q3D fitting procedure (§ 4.1) and the consequent single blade Q3D optimization (§ 4.2). Conventional Hub, Mid and Tip sections will be considered. Starting geometries come out from 1D meanline solution, as they are automatically generated without any previous PD optimization.

Since Q3D analyses represent almost all the time of the optimization procedures a notable time saving is reached introducing parallelization (when possible) in the optimization routines, so that multiple Q3D analyses can be performed simultaneously. For all the simulations the same 8 core-based system with no bottlenecks in memory was used. Attention was also paid to load equally the machine in order not to give a bias to the obtained results.

### 4.1 Fitting

As anticipated in Part 2, the "fitting" is a fundamental step of the PD required to ensure a reliable geometry optimization [11]. Generally speaking, it consists in overlapping the isentropic Mach profiles that contour the blades provided by the 3D analyses to those obtained by the Q3D simulations, thus imposing the same flow field around the turbine section. This result is achieved by modifying five different variables which define the flow channel characteristics: the inlet flow angle, three channel diffusion factors, and a coefficient for the exit static pressure. Acting on these parameters, the fitting procedure aims at minimizing the isentropic Mach error, which is calculated using the standard deviation between the 3D and the Q3D Mach number values on both suction and pressure side.

It is apparent that the fitting is a challenging mono-objective optimization problem, characterized by jagged and very large boundaries. Furthermore, the solution may not be unique and the domain space usually presents many minima with close objective function values. Therefore, the choice of the correct optimization algorithm to use is difficult: some methods can be severely affected by solution convergence, while intrinsically multi-objective (GeDEA and ANN) ones are not well exploited.

In the present work, the fitting is carried out on the three sections (hub, mid, and tip) of the test-case blade. The algorithms presented in Chapter 3 were tested and compared in order to assess their performances with respect to the minimization problem under consideration. LOH+ANN method was not taken into account because the vastness and uncertainty of the design space make the LOH sampling ineffective, unless adopting a prohibitive number of points.

Due to convergence error, it happens that identical set of input parameters do not always provide the same distribution of the isentropic Mach number, which means that the convergence precision of the Q3D analysis affects the results of the numerical simulation. This fact can be interpreted as noise that affects more the methods based on differences from small perturbations (mainly the GD). The objective function was specifically created to include not only the Mach number error, but also the convergence error. Hence all the optimization algorithms prioritize the solutions whose results are more trustworthy.

### 4.1.1 Results and analysis

Figure 2 shows an example of the results obtained during the fitting of the hub section. It is easy to note that Q3D Mach contours (gold) are almost perfectly overlapped on 3D ones (blue), while the initial condition (green) is very far from the final outline as usual for hub sections due to secondary 3D effects. Obviously, geometries remain unchanged during all fitting operations.

In order to reduce the impact of the random component, all the analyses were repeated six times. Table 2 reports the comparison of the algorithms in terms of computational resources used. The total time was kept constant (one hour) for each random-based optimization method; however, the deterministic GD always reached a local minimum before the allowed time.

Figure 3 shows the trend of the objective function with respect to the time for all the optimization algorithms for hub and tip sections. Note that the plotted objective function (error ratio) is the error divided by a small reference value and the error ratio axis has been cut to improve readability. The trends of the objective function with respect to number of iterations, here not plotted, do not present remarkable differences from the ones depending on time in Figure 3. This was expected, thanks to a fine code parallelization. It is also worth noting that the plot reported in Figure 3 compare the trend of the objective function obtained by the IRW to those provided by three parallelized algorithms (GD, AsBeC and GeDEA).

The IRW is characterized by the poorest performance level, not only due to the lower number of analyses performed within the prescribed time. In fact, AsBeC and GeDEA are designed to better exploit the high number of simulations they can compute. It has been tested that even running 8 IRW instances in parallel the best result is still worse than both the average of the two evolutionary algorithms. AsBeC seems to have better exploring qualities than others, finally reaching lower objective values in hub (shown) and mid (not shown) cases. In addition, AsBeC final results are less affected by the random component as it can be deduced by looking at the mean standard deviation reported in Table 1.

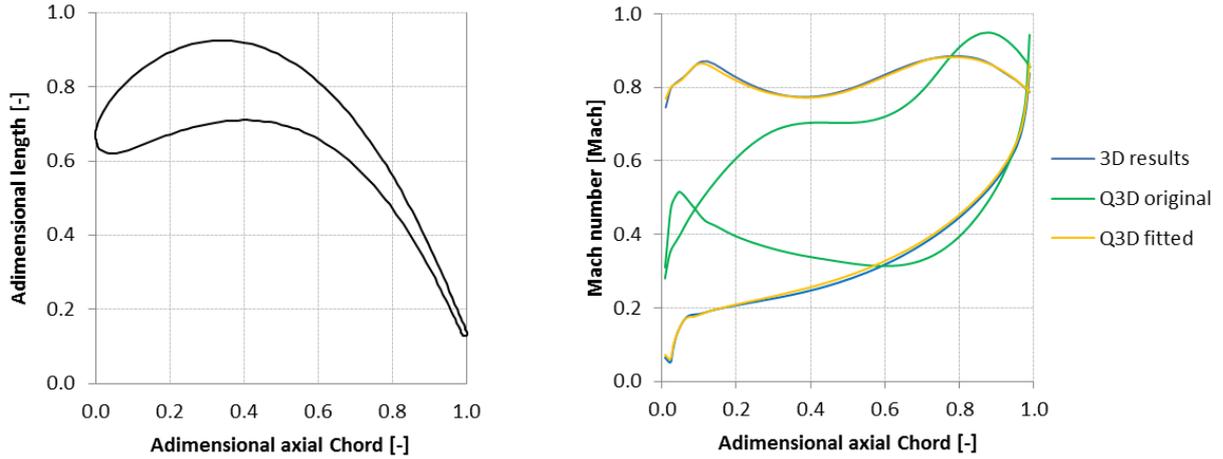

Figure 2: Hub blade section (on the left) and Mach profiles before and after fitting (on the right).

|  | GD | IRW | AsBeC | GeDEA |
|---|---|---|---|---|
| Total time (minutes) | 6 (mean) | 60 | 60 | 60 |
| Mean total # of Q3D analyses performed | 110 | 193 | 1064 | 1104 |
| # of threads used per block | 5-6 | 1 | 8 | 8 |
| Mean # blocks of parallel analyses | 20 | 193 | 133 | 138 |
| Mean final standard deviation on error ratio | - | 0.1 | 0.007 | 0.04 |

Table 1: Comparison of mean resources used among the algorithms for the fitting of the 3 sections

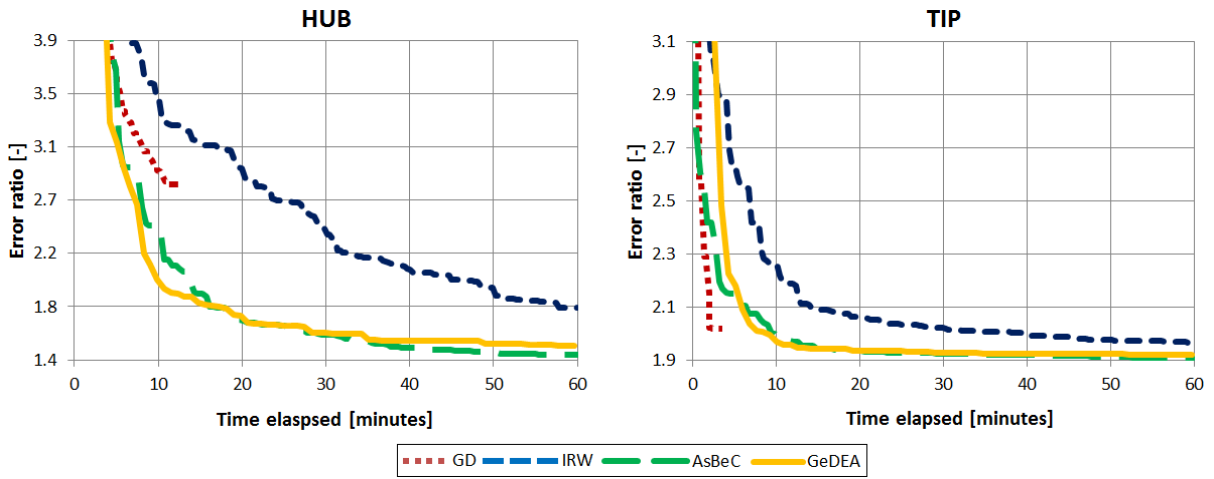

Figure 3: Performance comparisons of the algorithms in hub and tip sections

## 4.2 Optimization

The geometric optimization of LPT airfoils is the core of the PD phase. Focusing on the specific study, the parameters domain is six-dimensional and the variables are namely: the tangential chord, the unguided turning, the inlet blade angle, the inlet wedge angle, the leading edge radius and its eccentricity. The optimization of the high speed blade airfoils is multidisciplinary and essentially aero-mechanical. Three different objective functions were properly defined so that they have to be all maximized. In particular, the first objective is represented by the aerodynamic efficiency (η), the second one is connected to section area (Area Obj.) in order to fulfil the structural requirements; the third function (M+C Obj.) combines together the convergence error (to prioritize trustworthy solutions) and an objective on the Mach numbers along the profile suction side, that have to be limited for aerodynamic reasons. In the case of the GD, the IRW and the AsBeC, these three objectives were gathered together in only one function by means of a weighted sum.

In this multi-objective framework the Pareto front (possible with ANN and GeDEA) can be very useful, allowing to give more importance to some objective functions without re-performing the whole optimization.

The five different optimization methods are compared as done in the fitting case, performing six different instances for each section. In addition, the Pareto Front is presented and analysed from the physical point of view.

*4.2.1 Results and analysis*

This time all the analyses were stopped after 8 hours (a working day), with the only exception of the GD which reaches its final solution after 28 minutes on average, as reported in Table 2. In Figure 4, the objective functions for the hub and the mid are reported with respect to the time. As far as the multi-objective optimization algorithms (GeDEA and ANN) are concerned, Figure 4 plots the trend of the best values in terms of weighted sum of the three objective functions. Thanks to this expedient, it is possible to appreciate the performances of all the algorithm used, even though the comparison is not completely consistent. The right side of Figure 4 reports the airfoils and isentropic Mach trends of the optimized configuration (black solid lines) with respect to the baseline (red solid lines).

As for the fitting, the IRW method is the slowest due to the same reasons explained in §4.1. The AsBeC algorithm combines global and local search abilities in a very efficient way, being faster than IRW and overcoming even the GeDEA after some time. The LOH+ANN scheme turns out to be the best overall despite the long training time required. It also runs less Q3D simulations than the AsBeC and the GeDEA, as apparent from Table 2. On the other hand, the GeDEA seems to be unable at finding better solutions after some hours of evolutionary processes, probably because of elitism. All the previous observations are also valid for the tip section, whose results are omitted for the sake of compactness.

The remarkable performance of the ANN is mostly linked to the problem under study, in which the initial configuration was very different from the optimized one. In such a situation, the accurate LOH sampling is fundamental for the investigation of the feasible domain. In order to discuss an antithetical example, a midspan section closer to the optimum one was optimized apart for 60 minutes using all the algorithms. The results are shown in Figure 5: compared to Figure 4 the differences are relevant. In this case, the LOH+ANN performance is evidently the worst, since the LOH is totally ineffective when sampling in a large domain around an already good initial solution. Following the same wide-exploration principle, also GeDEA degrades its quality. On the contrary, GD and AsBeC are built to focus on the local search and then proves to be the most beneficial.

However, LOH+ANN and GeDEA provide the Pareto Fronts which can be very useful to ensure a quick and excellent design. A practical illustration is offered in Figure 6, in which both the whole final GeDEA Pareto Front (greyscale for M+C Obj.) and an intermediate step of the LOH+ANN Pareto front (yellow-to-blue scale for M+C Obj., zoomed inside the picture) are shown. It is worth remembering that both the algorithms thicken their Pareto near points that have higher weighted sum of objectives. Moreover, the diversity within the population in the GeDEA helps to obtain widespread Pareto Front. On the contrary, LOH+ANN Paretos narrow the more the refinement is advanced, but the ANN metamodel allows to enrich the specialized Pareto with much more points. Hub section is selected as the key example since it has to balance mechanical and aerodynamic goals. In general, Area Obj. results to be inversely proportional to efficiency while M+C Obj. seems positively correlated with efficiency.

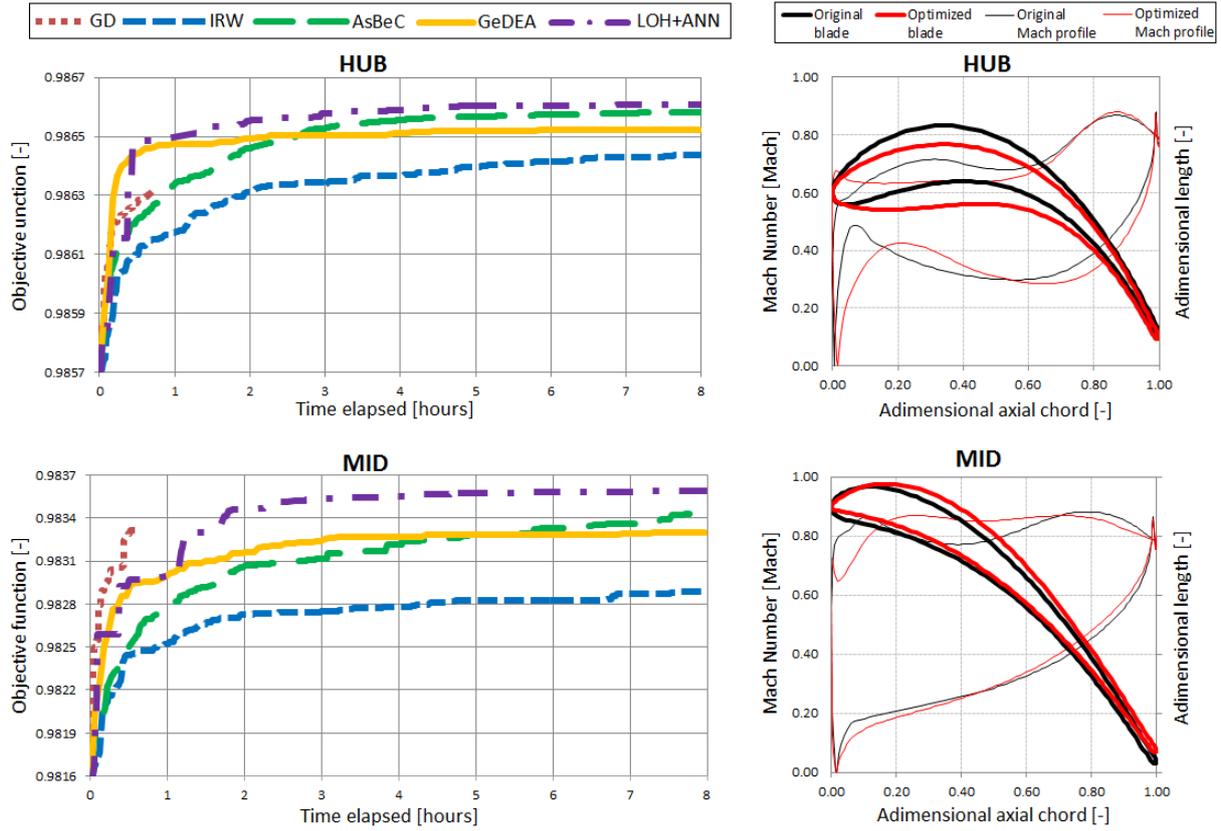

Figure 4: Algorithms performance comparisons and final optimized blade in different sections

|  | GD | IRW | AsBeC | GeDEA | ANN-LOH |
|---|---|---|---|---|---|
| Total time (minutes) | 28 (mean) | 480 | 480 | 480 | 480 |
| Mean total # of Q3D analyses performed | 203 | 532 | 3232 | 3360 | 2816 |
| # of threads used per block | 6-8 | 1 | 8 | 8 | 8 |
| Mean # blocks of parallel analyses | 29 | 532 | 404 | 420 | 352 |
| Mean final standard deviation on objective | - | $1 \cdot 10^{-4}$ | $6 \cdot 10^{-5}$ | $4 \cdot 10^{-5}$ | $3 \cdot 10^{-5}$ |

Table 2: Comparison of resources used among the algorithms for the optimization of the 3 sections

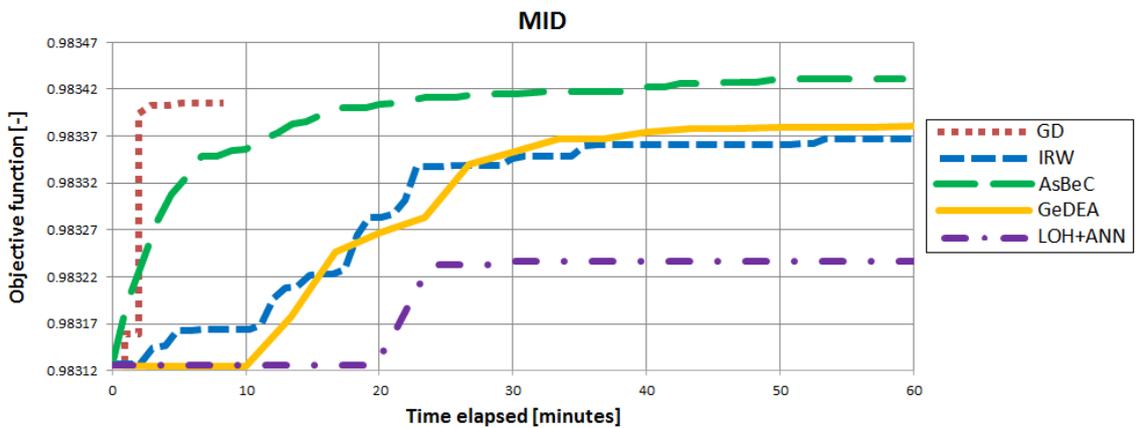

Figure 5: Performance comparisons of the algorithms starting from an already optimized blade

Three different solutions belonging to the ANN Pareto front (named Case 1, Case 2 and Case 3), have been tested and compared (right side of Figure 6). Specifically, the Case 1 (blue) is an high efficiency configuration, that limits the diffusive phenomena, lowers the Mach number and promotes the flow acceleration; Case 2 (green) is the best configuration with balanced weights (same of previous optimizations), and an intermediate Mach profile between Case 1 and Case 3; Case 3 (gold) is a high resistance configuration, that increases the area at the expense of higher Mach number and strong diffusion, mainly on the pressure side.

The Pareto exploitation reveals how much easier and faster can be the blade design moving among optimal configurations without re-performing any further optimization. For the specific example of hub section, Case 3 is recommended because its lower η is largely dominated by its structural strength that allows larger safety margin on breakages.

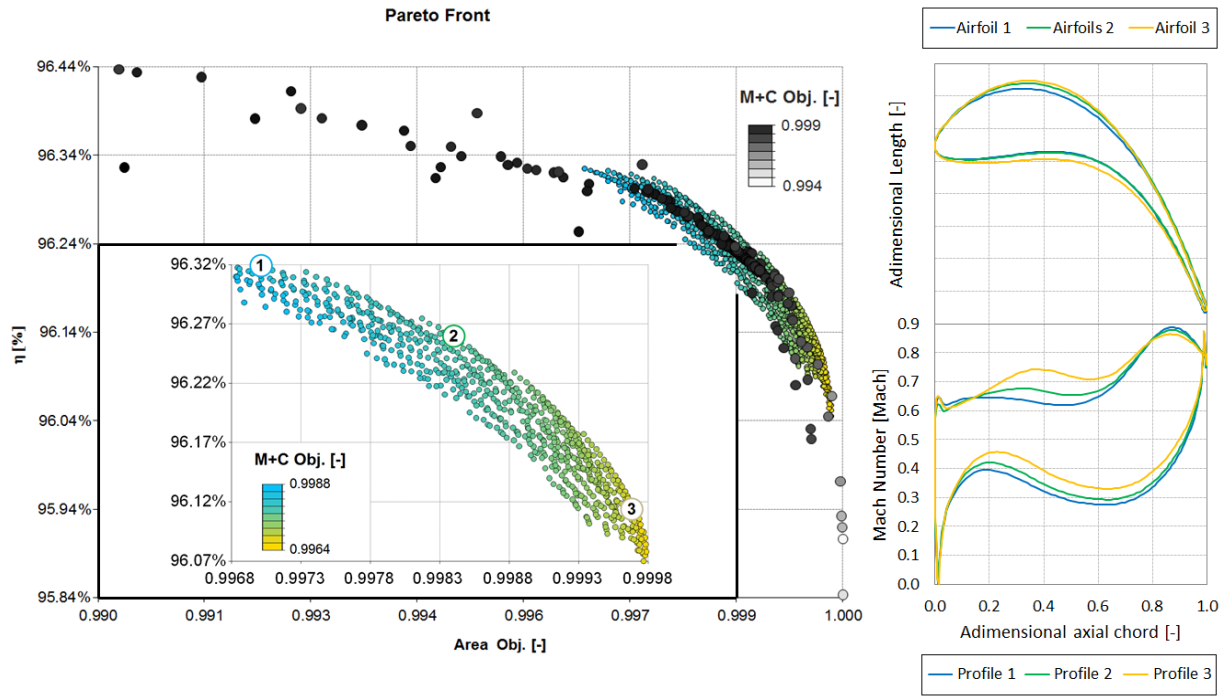

Figure 6: 3D Pareto Fronts for GeDEA and ANN and 3 optimal configurations tested (hub section)

## 5     CONCLUDING REMARKS

This paper concerned the application and comparison of different optimization methods for the PD design of a real-like and peculiar LPT rotor blade. Quasi mono-objective optimization with large parameter space (fitting operation) was addressed as well as typical multi-objective one, with main aero-mechanical goals but limited feasibility domain (geometric optimization). In the first case, both evolutionary algorithms turn out to be the best, but the GeDEA needs time to accurately set the range of the boundaries, while the AsBeC shows a better local search attitude, with the possibility of finding lower minima. In the second case, the wide-exploring algorithms are the best choice among the others, especially when the optimization starts from a poor initial configuration. Nevertheless, LOH+ANN is not suitable for high dimensionality in parameters and GeDEA tends to freeze the genetic of the individuals due to elitism when advancing in time. Furthermore, optimizing an already good solution leads to worsen both of them (especially LOH+ANN) while AsBeC ad GD works very well. However, LOH+ANN and GeDEA should be privileged for multi-objective and multi-disciplinary problems since they provide the Pareto Front.

## 6     ACKNOWLEDGEMENTS

The authors would like to thank M. Marconcini and M. Giovannini at the University of Florence for the support on the TRAF code for CFD analyses.